\documentclass[12pt]{article}
 \usepackage{latexsym}
 \usepackage{amssymb}
 \usepackage{graphicx}
\usepackage{color}
 \newtheorem{Theorem}{Theorem}

\newtheorem{Proposition}{Proposition}

\newtheorem{Corollary}{Corollary}
\newtheorem{Remark}{Remark}
\newtheorem{Example}{Example}

\newenvironment{example}[1][Example]{\begin{trivlist}
\item[\hskip \labelsep {\bfseries #1}]}{\end{trivlist}}

\newcommand{\A}{{\cal A}}
\newcommand{\B}{{\cal B}}
\newcommand{\C}{{\cal C}}

\newcommand{\HH}{{\cal H}}

\newcommand{\OO}{{\cal O}}
\newcommand{\PP}{{\cal P}}
\newcommand{\X}{{\cal X}}

\newcommand{\uu}{{\bf u}}
\newcommand{\vv}{{\bf v}}
\newcommand{\x}{{\bf x}}
\newcommand{\y}{{\bf y}}

\newcommand{\e}{{\bf e}}
\newcommand{\0}{{\bf 0}}

\newcommand{\qed}{\nobreak \ifvmode \relax \else
      \ifdim\lastskip<1.5em \hskip-\lastskip
      \hskip1.5em plus0em minus0.5em \fi \nobreak
      \vrule height0.75em width0.5em depth0.25em\fi}

\def \ep{\hbox{ }\hfill$\Box$}

\begin{document}
\title{SOS-Hankel Tensors: Theory and Application}

\author{Guoyin Li
\footnote{Department of Applied Mathematics, University of New South
Wales, Sydney 2052, Australia. E-mail: g.li@unsw.edu.au (G. Li).
This author's work was partially supported by Australian Research
Council.} \quad Liqun Qi\footnote{Department of Applied
 Mathematics, The Hong Kong Polytechnic University, Hung Hom, Kowloon, Hong
 Kong.
 E-mail: maqilq@polyu.edu.hk (L. Qi).  This author's work was partially supported by the Hong Kong Research Grant Council (Grant No. PolyU
502510, 502111, 501212 and 501913). } \quad Yi Xu \footnote{Department of Mathematics, Southeast
University, Nanjing, 210096 China.   Email: yi.xu1983@gmail.com (Y. Xu).}}

\date{\today} \maketitle

\begin{abstract}
\noindent  
Hankel tensors arise from signal processing and some other applications.
SOS (sum-of-squares) tensors are positive semi-definite symmetric tensors, but not vice versa.  The problem for determining an even order symmetric tensor is an SOS tensor or not is equivalent to solving a semi-infinite linear programming problem, which can be done in polynomial time.  On the other hand, the problem  for determining an even order symmetric tensor is positive semi-definite or not is NP-hard.   In this paper, we study SOS-Hankel tensors.   Currently, there are two known positive semi-definite Hankel tensor classes: even order complete Hankel tensors and even order strong Hankel tensors.   We show complete Hankel tensors are strong Hankel tensors, and even order strong Hankel tensors are SOS-Hankel tensors.    We give several examples of positive semi-definite Hankel tensors, which are not strong Hankel tensors.   However, all of them are still SOS-Hankel tensors.   Does there exist a positive semi-definite non-SOS-Hankel tensor?    The answer to this question remains open.  If the answer to this question is no, then the problem for determining an even order Hankel tensor is positive semi-definite or not is solvable in polynomial-time.     An application of SOS-Hankel tensors to the positive semi-definite tensor completion problem is discussed. We present an ADMM algorithm for solving this problem.  Some preliminary numerical results on this algorithm are reported.

\noindent {\bf Key words:}\hspace{2mm} Hankel tensors, generating vectors, sum of squares, positive semi-definiteness,
generating functions, tensor completion. \vspace{3mm}

\noindent {\bf AMS subject classifications (2010):}\hspace{2mm}
15A18; 15A69
  \vspace{3mm}

\end{abstract}

\section{Introduction}
\hspace{4mm}

In general, most tensor problems are NP-hard \cite{HL}.  However, most tensor problems in applications have special structures, and they are not NP-hard.   In the last few years, there are a number of papers on the Perron-Frobenius theory for nonnegative tensors, and algorithms for computing the largest eigenvalue of a nonnegative tensor \cite{CPZ, CQZ, NQZ}.  In particular, in \cite{HLQS}, it was shown that the problem for computing the largest H-eigenvalue of an essentially nonnegative tensor, which includes the problem for computing the largest eigenvalue of a nonnegative tensor as a special problem, can be solved by solving a semi-definite linear programming problem.  Thus, this problem is polynomial-time solvable, not NP-hard.    This method can be used to find the smallest H-eigenvalue of a Z tensor, thus to be used to determine a given Z tensor is an M tensor or not, while an even order symmetric M tensor is positive semi-definite \cite{ZQZ}.

The problem for determining a given even order symmetric tensor is positive semi-definite or not has important applications in engineering and science \cite{Hi, Qi, QYW, QYX, Re, Sh}.   In general, this problem is NP-hard.  However, for special tensor classes, such as even order symmetric Z tensors,
as pointed above, this problem is polynomial time solvable.
In 2014, more classes of structured tensors have been identified, either such
tensors are easily to be identified, and they are positive semi-definite or positive definite in the even order symmetric case, or there are easily checkable conditions to identify such tensors are positive semi-definite or not.  The former includes Hilbert tensors \cite{SQ}, diagonally dominated tensors \cite{QS}, B tensors \cite{QS}, double B tensors \cite{LL}, quasi-double B tensors \cite{LL} and H$^+$ tensors \cite{LWZZL}.  The latter includes Cauchy tensors \cite{CQ}.

In \cite{HLQ}, a new class of positive semi-definite tensors, called SOS (sum-of-squares) tensors, was introduced.  SOS tensors are positive semi-definite tensors, but not vice versa.   SOS tensors are connected with SOS polynomials, which are significant in polynomial theory \cite{Hi, Re, Sh} and polynomial optimization \cite{Las, Lau}.   In particular, as stated above, the problem to identify a given general even order symmetric tensor is positive semi-definite or not is NP-hard, while the problem to identify a given general even order symmetric tensor is SOS or not is equivalent to solving a semi-definite linear programming problem \cite{Las, Lau}, thus not NP-hard, but polynomial time solvable.   However, no special structured tensor class was identified to be SOS tensors in \cite{HLQ}.

Hankel tensors arise from signal processing and some other applications \cite{BB, DQW, PDV, Qi15}.   Hankel tensors are symmetric tensors.  In \cite{Qi15}, positive semi-definite Hankel tensors were studied.   Each Hankel tensor is associated with an Hankel matrix. If that Hankel matrix is positive semi-definite, then the Hankel tensor is called a strong Hankel tensor.   It was proved that an even order strong Hankel tensor is positive semi-definite.
A symmetric tensor is a Hankel tensor if and only if it has a Vandermonde decomposition.   If the coefficients of that Vandermonde decomposition are nonnegative, then the Hankel tensor is called a complete Hankel tensor.    It was proved that an even order complete Hankel tensor is also positive semi-definite.  The relation between strong and complete Hankel tensors was not known in \cite{Qi15}.  An example of a positive semi-definite Hankel tensor, which is neither strong nor complete Hankel tensor was also given in \cite{Qi15}.

In this paper, we study positive semi-definite and SOS Hankel tensors.   We introduce completely decomposable tensors.  Even order completely decomposable tensors are SOS tensors, but not vice versa.    We show that complete Hankel tensors are strong Hankel tensors, while strong Hankel tensors are completely decomposable.  Thus, both even order complete Hankel tensors and even order strong Hankel tensors are SOS Hankel tensors.   From now on, we call SOS Hankel tensors as SOS-Hankel tensors.  Then we show that for any even order $m=2k\ge 4$, there are SOS-Hankel tensors, which are not completely decomposable.  This includes the example given in \cite{Qi15}.  We also found some other examples of SOS-Hankel tensors which are not strong Hankel tensors.  Does there exist a positive semi-definite non-SOS Hankel tensor?    The answer to this question remains open.   If the answer to this question is no, then the problem for determining an even order Hankel tensor is positive semi-definite or not is solvable in polynomial-time.

The remainder of this paper is organized as follows.  Some preliminary knowledge is given in the next section.   In Section 3, we introduce completely decomposable tensors and discuss their properties.   We prove completely Hankel tensors are strong Hankel tensors, and strong Hankel tensors are completely decomposable in Section 4.   Some SOS-Hankel tensors which are not strong Hankel tensors are given in Section 5.  Then we raise the question on positive semi-definite non-SOS Hankel tensors and have some discussion.   Finally, in Section 6, we discuss an application of SOS-Hankel tensors to the positive semi-definite tensor completion problem.  We present an ADMM algorithm for solving this problem.  Some preliminary numerical results on this algorithm are reported.

\section{Preliminaries}
\hspace{4mm}

Throughout this paper, we use small letters $a, b, \alpha, \beta, \lambda, \cdots$, for scalars, bold small letters $\x, \y, \uu, \vv, \cdots$, for vectors, capital letters $A, B, C, I, \cdots$, for matrices, calligraphic letters $\A, \B, \C, \HH, \cdots$, for tensors.  We use $\0$ to denote the zero vector in $\Re^n$, and $\e_i$ to denote the $i$th unit vector in $\Re^n$, for $i \in [n]$.   We use $x_i$ to denote the $i$th component of $\x$.

Denote that $[n] := \{ 1, \cdots, n \}$.   Throughout this paper, $m, n$ and $k$ are integers, and $m, n \ge 2$.   A tensor $\A = (a_{i_1\cdots i_m})$ of order $m$ and dimension $n$ has entries $a_{i_1\cdots i_m}$ with $i_j \in [n]$ for $j \in [m]$.  Tensor $\A$ is said to be a symmetric tensor if its entries
$a_{i_1\cdots i_m}$ is invariant under any index permutation.   Denote the set of all the real symmetric tensors of order $m$ and dimension $n$ by ${\rm S}_{m, n}$.  Then ${\rm S}_{m, n}$ is a linear space.  Throughout this paper, we only discuss real symmetric tensors.  We use $\| \A \|$ to denote the Frobenius norm of tensor $\A = (a_{i_1\cdots i_m})$, i.e., $\| \A \| = \sum_{i_1, \cdots, i_m \in [n]} a_{i_1 \cdots i_m}^2$.

Let $\x \in \Re^n$.   Then $\x^{\otimes m}$ is a rank-one symmetric tensor with entries $x_{i_1}\cdots x_{i_m}$.  For $\A \in {\rm S}_{m, n}$ and $\x \in \Re^n$, we have a homogeneous polynomial $f(\x)$ of $n$ variables and degree $m$,
\begin{equation} \label{e2.0}
f(\x) = \A \x^{\otimes m} \equiv \sum_{i_1, \cdots, i_m \in [n]}a_{i_1\cdots i_m}x_{i_1}\cdots x_{i_m}.
\end{equation}
Note that there is a one to one relation between homogeneous polynomials and symmetric tensors.

If $f(\x) \ge 0$ for all $\x \in \Re^n$, then homogeneous polynomial $f(\x)$ and symmetric tensor $\A$ are called {\bf positive semi-definite} (PSD).   If $f(\x) > 0$ for all $\x \in \Re^n, \x \not = \0$, then $f(\x)$ and $\A$ are called {\bf positive definite}.  Clearly, if $m$ is odd, there is no positive definite symmetric tensor and there is only one positive semi-definite tensor $\OO$.  Thus, we assume that $m = 2k$ when {  we discuss positive definite and semi-definite tensors (polynomials). Positive semi-definite polynomials, or called nonnegative polynomials,} have important significance in mathematics, and are connected with the 17th Hilbert problem \cite{Hi, Re, Sh}.   The concepts of positive semi-definite and positive definite symmetric tensors were introduced in \cite{Qi}.  They have wide applications in science and engineering \cite{Qi, QY, QYW, QYX}. Denote the set of all positive semi-definite symmetric tensors in ${\rm S}_{m, n}$ by ${\rm PSD}_{m, n}$.  Then it is a closed convex cone and has been studied in \cite{HLQ, QY}.

Let $m=2k$.   If $f(\x)$ can be decomposed to the sum of squares of polynomials of degree $k$, then $f(\x)$ is called a {\bf sum-of-squares (SOS) polynomial}, and the corresponding symmetric tensor $\A$ is called an {\bf SOS tensor} \cite{HLQ}.    SOS polynomials play a central role in the modern theory of polynomial optimization \cite{Las, Lau}.   SOS tensors were introduced in \cite{HLQ}. Clearly, an SOS polynomial (tensor) is a PSD polynomial, but not vice versa.   Actually, this was shown by young Hilbert \cite{CL, Hi, Mo, Re, Sh} that for homogeneous polynomial, only in the following three cases, a PSD polynomial definitely is an SOS polynomial: 1) $n = 2$; 2) $m = 2$; 3) $m=4$ and $n=3$.
For tensors, the second case corresponds to the symmetric matrices, i.e., a PSD symmetric matrix is always an SOS matrix.   Hilbert proved that in all the other possible combinations of $m=2k$ and $n$, there are non-SOS PSD homogeneous polynomials \cite{Sh}.  The most well-known non-SOS PSD homogeneous polynomial is the Motzkin polynomial \cite{Mo}
$$f_M(\x) = x_3^6 + x_1^2x_2^4 + x_1^4x_2^2 - 3 x_1^2x_2^2x_3^2.$$
By the Arithmetic-Geometric inequality, we see that it is a PSD polynomial.  But it is not an SOS polynomial \cite{Re, Sh}.   The other two non-SOS PSD homogeneous polynomials with small $m$ and $n$ are given by Choi and Lam \cite{CL}
$$f_{CL1}(\x) = x_4^4 + x_1^2x_2^2 + x_1^2x_3^2 + x_2^2x_3^2 - 4x_1x_2x_3x_4$$
and
$$f_{CL2}(\x) = x_1^4x_2^2 + x_2^4x_3^2 + x_3^4x_1^2 - 3x_1^2x_2^2x_3^2.$$
Denote the set of all SOS tensors in ${\rm S}_{m, n}$ by ${\rm SOS}_{m, n}$.  Then it is also a closed convex cone \cite{HLQ}.   Thus, ${\rm SOS}_{m, 2} = {\rm PSD}_{m, 2}$, ${\rm SOS}_{2, n} = {\rm PSD}_{2, n}$ and ${\rm SOS}_{4, 3} = {\rm PSD}_{4, 3}$.   But for other $m =2k \ge 4, n \ge 3$, we have ${\rm SOS}_{m, n} \subsetneq {\rm PSD}_{m, n}$.

{  Let $\vv = (v_0, \cdots, v_{(n-1)m})^\top$.   Define $\A = (a_{i_1\cdots i_m}) \in {\rm S}_{m, n}$ by
\begin{equation} \label{e2.1}
a_{i_1\cdots i_m} = v_{i_1+\cdots + i_m-m},
\end{equation}
for $i_1, \cdots, i_m \in [n]$.}   Then $\A$ is a {\bf Hankel tensor} and $\vv$ is called the {\bf generating vector} of $\A$.   In (\ref{e2.0}), if $\A$ is a Hankel tensor, {  then} $f(\x )$ is called a {\bf Hankel polynomial}.   We see that  a sufficient and necessary condition for $\A = (a_{i_1\cdots i_m}) \in S_{m, n}$ to be a Hankel tensor is that whenever $i_1+\cdots + i_m = j_1+ \cdots + j_m$,
\begin{equation} \label{e2.2}
a_{i_1\cdots i_m} = a_{j_1\cdots j_m}.
\end{equation}
Hankel tensors arise from signal processing and other applications \cite{BB, DQW, PDV, Qi15}.

By (\ref{e2.2}), the three non-SOS PSD polynomials $f_M(\x )$, $f_{CL1}(\x )$  and $f_{CL2}(\x )$ are not Hankel polynomials.   These three polynomials are still non-SOS PSD polynomials if we switch the indices of their variables.   However, if we switch the indices of their variables and add some terms to them to make them become Hankel polynomials, then the resulted Hankel polynomials are not positive semi-definite.   Hence, no non-SOS PSD Hankel polynomials are resulted.    
There are other examples of non-SOS PSD polynomials \cite{Re}.   None of them are Hankel polynomials.


Note that for $n = 2$, all symmetric tensors are Hankel tensors.

\section{Completely Decomposable Tensors}
\hspace{4mm}

Let $\A \in {\rm S}_{m, n}$.  If there are vectors $\x_j \in \Re^n$ for $j \in [r]$ such that
$$\A = \sum_{j \in [r]} \x_j^{\otimes m},$$
then we say that $\A$ is a {\bf completely $r$-decomposable tensor}, or a {\bf completely decomposable tensor}.   If $\x_j \in \Re^n_+$ for all $j \in [r]$, then $\A$ is called a {\bf completely positive tensor} \cite{QXX}.  If $\x_j \in \Re^n_+$ for all $j \in [r]$, then $\A$ is called a {\bf completely positive tensor}. {
We note that } any odd order symmetric tensor is completely decomposable.  Hence, it does not make sense to discuss odd order completely decomposable tensor.  But it is still meaningful to show an odd order symmetric tensor is completely $r$-decomposable, as this means that
its symmetric rank \cite{CGLM} is not greater than $r$.

Denote the set of all completely $r$-decomposable tensors in ${\rm S}_{m, n}$ by ${\rm CD}^r_{m, n}$, and the set of all completely decomposable tensors in ${\rm S}_{m, n}$ by ${\rm CD}_{m, n}$.  The following theorem will be useful in the next section.

\begin{Theorem}{\bf (Closedness and convexity of the tensor cones)} \label{t3.1}
 The following statements hold:
 \begin{itemize}
 \item[{\rm (i)}] If the order $m$ is an even number, then ${\rm SOS}_{m,n}$ and ${\rm PSD}_{m,n}$ are closed and convex cones.
 \item[{\rm (ii)}] If the order $m$ is an even number, ${\rm CD}_{m,n}^r$ is a closed cone for each integer $r$.
 \item[{\rm (iii)}] If the order $m$ is an odd number, then
 \begin{itemize}
  \item[{\rm (a)}] ${\rm CD}_{m,n}^{r}$ is either a nonconvex cone such that ${\rm CD}_{m,n}^{r}=-{\rm CD}_{m,n}^{r}$ or the whole symmetric tensor space ${\rm S}_{m,n}$. In particular, for any
 $r<n$, ${\rm CD}_{m,n}^{r}$ is a nonconvex cone such that ${\rm CD}_{m,n}^{r}=-{\rm CD}_{m,n}^{r}$;
  \item[{\rm (b)}] ${\rm CD}_{m,n}^{r}$ can be non-closed for some integer $r$.
 \end{itemize}
 \item[{\rm (iv)}] If the order $m$ is an even number, then ${\rm CD}_{m,n}$ is a closed and convex cone. If $m$ is an odd number, then ${\rm CD}_{m,n}={\rm S}_{m,n}$.
 \end{itemize}
 \end{Theorem}

 \noindent {\bf Proof.}

[Proof of {\rm (i)}] The closed and convexity of positive
semidefinite cone ${\rm PSD}_{m,n}$ can be directly verified from
the definition. The convexity of the SOS tensor cone ${\rm
SOS}_{m,n}$ also directly follows from the definition. To see the
closedness of the SOS tensor cone ${\rm SOS}_{m,n}$, let
$\A_k \in {\rm SOS}_{m,n}$ with $\mathcal{A}_k \rightarrow
\A \in {\rm S}_{m,n}$. So, $f_k(\x )=\A_k\x^{\otimes m}$ is
an SOS polynomial and $f_k(x) \rightarrow f(x)$ with
$f(\x )=\A \x^{\otimes m}$. Note from \cite[Corollary
3.50]{Lau} that the set of all sum-of-squares polynomials on $\Re^n$
(possibly nonhomogeneous) with degree at most $m$ is a closed cone. So, $f$ is also an SOS polynomial. Therefore, $\A$ is an SOS
tensor.

[Proof of {\rm (ii)}] To see the closedness of ${\rm CD}^r_{m,n}$,
we let $\A_k \in {\rm CD}^r_{m,n}$ with $\A_k
\rightarrow \A$. Then, for each integer $k$,  there
exist $\x^{k,j} \in \Re^n$, $j \in [r]$, such that
$\A_k = \sum_{j\in [r]} \left(\x^{k,j}\right)^{\otimes m}$. As
$\A_k \rightarrow \A$, $\left\{\|\A_k\|\right\}_{k
\in \mathbb{N}}$ is a bounded sequence. Note that
$\|\A_k\|^2 \ge    \sum_{j \in [r]}\left(\sum_{i=1}^n
\left(x^{k,j}_i\right)^m\right)^2$. So, $\{\x^{k,j}\}_{k \in \mathbb{N}}$,
$j=1,\ldots,r$, are bounded sequences. By passing to subsequences,
we can assume that $\x^{k,j} \ \rightarrow \x^j$, $j=1,\ldots,r$.
Passing to the limit, we have
\[
 \A = \sum_{j\in [r]} \left(\x^j\right)^{\otimes m} \in {\rm CD}^r_{m,n}.
\]
Thus, the conclusion follows.

[Proof of {\rm (iii) (a)}] First of all, it is easy to see that ${\rm CD}_{m,n}^r$ is a cone.  From the definition of ${\rm CD}_{m,n}^r$, if $\A \in {\rm CD}_{m,n}^r$, then {  there exist  $\uu_1,\ldots,\uu_r \in \Re^n$} such that
$\A =\sum_{i=1}^r u_i^{\otimes m}$.
As $m$ is odd, $-\mathcal{A}=\sum_{i=1}^r (-u_i)^{\otimes m}$, and so, $-\mathcal{A} \in {\rm CD}_{m,n}$. So, ${\rm CD}_{m,n}^r = -{\rm CD}_{m,n}^r$.  To finish the proof, it suffices to show that if $CD_{m,n}^r$ is convex, then it must be
the whole space ${\rm S}_{m,n}$. To see this, we suppose that ${\rm CD}_{m,n}^r$ is convex. Then, ${\rm CD}_{m,n}^r$ is a convex cone with ${\rm CD}_{m,n}^r = -{\rm CD}_{m,n}^r$. This implies that ${\rm CD}_{m,n}^r$ is a subspace. Now, take $\A \in ({\rm CD}_{m,n}^r)^{\bot}$ where $L^{\bot}$ denotes the orthogonal complement in ${\rm S}_{m,n}$ of the subspace $L$. Then, $\langle \A, \x^{\otimes m} \rangle=\A \x^{\otimes m}=0$ for all $\x \in \Re^n$. This shows that $\A=\OO$. So, $({\rm CD}_{m,n}^r)^{\bot}=\{\OO \}$, and hence ${\rm CD}_{m,n}^r={\rm S}_{m,n}$.

Now, let us consider the case when $r<n$. To see that assertion that ${\rm CD}_{m,n}^{r}$ is a nonconvex cone such that ${\rm CD}_{m,n}^{r}=-{\rm CD}_{m,n}^{r}$, we only need to show ${\rm CD}_{m,n}^r \neq {\rm S}_{m,n}$. To see this, let $\e_i$ denotes the vector whose $i$th coordinate is one and the others are zero, $i \in [n]$, and define
$\overline{\A}:=\sum_{i\in [n]} \e_i^{\otimes m}$. In particular, we see that $\overline{\A} \e_j^{\otimes (m-1)}=\e_j$ for all $j\in [n]$.
We now show that $\overline{\A} \notin {\rm CD}_{m,n}^r$. Otherwise, {  there exist  $\uu_1,\ldots,\uu_r \in \Re^n$} such that
$\overline{\A}=\sum_{i\in [r]} (\uu_i)^{\otimes m}$.
This implies that, for each $j \in [n]$,
\[
\e_j=\overline{\A} \e_j^{\otimes (m-1)}=\sum_{i\in [r]} \alpha_{ij} \uu_i \mbox{ with } \alpha_{ij}=\left((\uu_i)^\top \e_j\right)^{m-1} \in \Re.
\]
This implies that $\e_j \in {\rm span}\{\uu_1,\ldots,\uu_r\}$ for all $j \in [n]$, and so, ${\rm dim}{\rm span}\{\uu_1,\ldots,\uu_r\} \ge n$. This is impossible as we assume that $r<n.$

[Proof of {\rm (iii) (b)}] We borrow an example in \cite{CGLM} to
illustrate the possible non-closedness of the cone ${\rm
CD_{m,n}^r}$ in the case where $m$ is odd \footnote{This example was
used to show the set consisting of all tensors with symmetric rank
over the complex field less or equal to a fixed number $l$ can be
non-closed. Here we show that this example can also be used to show
the non-closedness of complete $r$-decomposable tensor cones with
odd order.}.  Let $m=3$ and $n=2$. For each $\epsilon>0$, let
\[
f_{\epsilon}(x_1,x_2)=\epsilon^2(x_1+\epsilon^{-1}x_2)^3+ \epsilon^2(x_1-\epsilon^{-1}x_2)^3
\]
and let $\A_{\epsilon}$ be the associated symmetric tensor, that is, $\A_{\epsilon}\x^{\otimes m}=f_{\epsilon}(\x )$ for all $\x \in \Re^2$. Clearly
$\A_{\epsilon} \in {\rm CD}_{3,2}^2$.
Then, $f_{\epsilon}\rightarrow f$ as $\epsilon \rightarrow 0$ where
$f(x_1,x_2)=6x_1x_2^2$.
Let $\A$ be the associated symmetric tensor of $f$. Then, we have $\A_{\epsilon} \rightarrow \A$. We now see that $\A \notin {\rm CD}_{3,2}^2$. To see this
we proceed by contradiction. Then there exist $a_1,b_1,a_2,b_2 \in \Re$ such that
\[
6x_1x_2^2=(a_1x_1+a_2x_2)^3+(b_1x_1+b_2x_2)^3.
\]
By comparing the coefficients, we have $a_1=-b_1$, $a_2=-b_2$ and $a_1a_1^2+b_1b_2^2=2$, which is impossible.

[Proof of {\rm (iv)}] From the definition, ${\rm CD}_{m,n}$ is a convex cone. As {  ${\rm CD}_{m,n} \subseteq {\rm S}_{m,n}$,  the dimension of ${\rm CD}_{m,n}$ is at most $I(m,n)$ where $I(m,n)=\left(\begin{array}{cc}
                                                                                                                                                                                               n+m-1 \\
                                                                                                                                                                                               n-1
                                                                                                                                                                                               \end{array}
\right)$}. By Carath\'{e}odory's theorem and the fact that ${\rm CD}_{m,n}$ is a convex cone, ${\rm CD}_{m,n}={\rm CD}_{m,n}^r$ with $r=I(m,n)+1$. Thus, {\rm (ii)} implies that ${\rm CD}_{m,n}$ is a closed convex cone if the order $m$ is an even number. We now consider the case when $m$ is an odd number.  Similar to the argument in {\rm (iii)}, we see that
${\rm CD}_{m,n}=-{\rm CD}_{m,n}$. So, ${\rm CD}_{m,n}$ is a subspace when $m$ is odd. Now, take $\A \in ({\rm CD}_{m,n})^{\bot}$. Then, $\A \x^{\otimes m}=0$ for all $\x \in \Re^n$. This shows that $\A=\OO$. So, $({\rm CD}_{m,n})^{\bot}=\{\OO \}$, and hence ${\rm CD}_{m,n}={\rm S}_{m,n}$.
 \ep

\begin{Remark}
The intricate non-closedness of the set which consists of all rank-$r$ symmetric tensors over the complex field was discovered and examined in detailed in \cite{CGLM}. This fact shows the
significant difference between a matrix and a tensor with order greater than $2$. Our results here shows that the same feature can happen for the complete $r$-decomposable symmetric tensors
where the order is an odd number.
\end{Remark}

As we have seen in Theorem \ref{t3.1}~{\rm (iii)(a)}, the cone ${\rm CD}_{m,n}^r$ is not convex when $m$ is an odd number and $r<n$. We now illustrate this by a simple example.
\begin{Example}{\bf (Illustrating the possible nonconvexity of ${\rm CD}_{m,n}^r$)}
Consider $m=3$, $n=2$ and $r=1$. Then, ${\rm CD}_{3,2}^1=\{\A \in {\rm S}_{3,2}: \A \x^{\otimes 3}=\left(\uu^\top \x\right)^3 \mbox{ for some } \uu \in \Re^2\}$. Let $\A_i \in {\rm S}_{3,2}$ be such that $\A_i \x^{\otimes 3}= x_i^3$, $i=1,2$. We now claim that $\frac{1}{2}\A_1+\frac{1}{2}\A_2 \notin {\rm CD}_{3,2}^1$. Otherwise, there exists $\uu=(u_1,u_2)^\top \in \Re^2$ such that
\[
\frac{1}{2}x_1^3+\frac{1}{2}x_2^3=(u_1x_1+u_2x_2)^3.
\]
Comparing the coefficients, we see that $u_1=u_2=\sqrt[3]{\frac{1}{2}}$ and $3u_1u_2^3=3u_1^3u_2=0$. This is impossible. So, we have $\frac{1}{2}\A_1+\frac{1}{2}\A_2 \notin {\rm CD}_{3,2}^1$, and hence ${\rm CD}_{3,2}^1$ is not convex.
\end{Example}

\section{Strong Hankel Tensors and Complete Hankel Tensors}
\hspace{4mm}

Suppose that $\A$ is  a
Hankel tensor defined by (\ref{e2.1}).    Let $A = (a_{ij})$ be an $\lceil
{(n-1)m+2 \over 2} \rceil  \times \lceil {(n-1)m+2 \over 2} \rceil$
matrix with $a_{ij} \equiv v_{i+j-2}$, where $v_{2\lceil {(n-1)m
\over 2} \rceil}$ is an additional number when $(n-1)m$ is odd. Then
$A$ is a Hankel matrix, associated with the Hankel tensor
$\A$.  Clearly, when $m$ is even, such an associated Hankel matrix is unique. Recall from \cite{Qi15} that $\A$ is called a {\bf strong Hankel tensor}
if {  there exists an associated Hankel matrix $A$ is positive semi-definite}. Thus,
whether a tensor is a strong Hankel tensor or not can be verified
by using tools from matrix analysis. It has also been shown in \cite{Qi15} that
$\A$ is a strong Hankel tensor if and only if it is a
Hankel tensor and there exists an absolutely integrable real valued
function $h:(-\infty, +\infty) \rightarrow [0,+\infty)$
 such that  its generating vector $\vv = (v_0, v_1, \ldots, v_{(n-1)m})^\top$ satisfies
\begin{equation}\label{strong_Hankel}
v_k=\int_{-\infty}^{\infty}t^k h(t) dt, \ \ \ k=0,1,\ldots,(n-1)m.
\end{equation}
Such a real valued function $h$ is called the {\bf generating function} of
the strong Hankel tensor $\A$.

 A vector $\uu = (1, \gamma, \gamma^2, \ldots, \gamma^{n-1})^\top$ for some $\gamma \in \Re$ is called a
Vandermonde vector \cite{Qi15}.  If tensor $\A$ has the form
\[
\A =\sum_{i\in [r]} \alpha_i(\uu_i)^{\otimes m},
\]
where $\uu_i$ for $i = 1, \ldots, r$, are all Vandermonde vectors,
then we say that $\mathcal{A}$ has a Vandermonde decomposition.  It
was shown in \cite{Qi15} that a symmetric tensor is a Hankel tensor
if and only if it has a Vandermonde decomposition.  If the
coefficients $\alpha_i$ for $i = 1, \ldots, r$, are all nonnegative,
then $\mathcal{A}$ is called a {\bf complete Hankel tensor} \cite{Qi15}.
Clearly, a complete Hankel tensor is a completely decomposable
tensor. Unlike strong Hankel tensors, there is no clear method to check whether a Hankel tensor is a
complete Hankel tensor or not,  as the
Vandermonde decompositions of a Hankel tensor are not unique.  It
was proved in \cite{Qi15} that even order strong or complete Hankel
tensors are positive semi-definite, but the relationship between
strong Hankel tensors and complete Hankel tensors was unknown there.

In the following, we discuss the relation between strong Hankel
tensors and complete Hankel tensors.

\begin{Proposition} \label{p4.1}
Let $m,n \in \mathbb{N}$.   All the strong Hankel tensors of order
$m$ and dimension $n$ form a closed convex cone.  All the complete
Hankel tensors of order $m$ and dimension $n$ form a convex cone. A
complete Hankel tensor is a strong Hankel tensor.   On the other
hand, whenever $m$ is a positive even number and $n \ge 2$, there is a strong Hankel tensor which is not a complete Hankel tensor.
\end{Proposition}

 \noindent {\bf Proof.} By definition, it is easy to see that all the complete
Hankel tensors of order $m$ and dimension $n$ form a convex cone.
Since each strong Hankel tensor of order $m$ and dimension $n$ is
associated with a positive semi-definite Hankel matrix and all such
positive semi-definite Hankel matrices form a closed convex cone,
all the strong Hankel tensors of order $m$ and dimension $n$ form a
closed convex cone.

Consider a rank-one complete Hankel tensor $\A =\uu^{\otimes
m}$, where $\uu = (1, \gamma, \gamma^2, \ldots, \gamma^{n-1})^\top$. Then
we see the generating vector of $\A$ is {  $\vv = (1, \gamma,
\gamma^2, \ldots, \gamma^{(n-1)m})^\top$}.  When $m$ is even, let $A$
be a Hankel matrix generated by $\vv$.   When $m$ is odd, let $A$ be a
Hankel matrix generated by $(1, \gamma, \gamma^2, \ldots,
\gamma^{(n-1)m})^\top$.   Then we see that $A$ is positive
semi-definite.   Thus, $\mathcal{A}$ is a strong Hankel tensor. From
the definition of complete Hankel tensors, we see that a complete
Hankel tensor is a linear combination of strong Hankel tensors with
nonnegative coefficients.  Since all the strong Hankel tensors of
order $m$ and dimension $n$ form a closed convex cone, this shows
that a complete Hankel tensor is a strong Hankel tensor.

Finally, assume that $m$ is a positive even number and $n \ge 2$.   Let $\A =\e_n^{\otimes m}$ where
$\e_n=(0,\ldots,0,1)\in \Re^n$.   It is easy to verify that
$\A$ is a strong Hankel tensor. We now see that
$\A$ is not a complete Hankel tensor. Assume on the
contrary. Then, there exist $r \in \mathbb{N}$, $\alpha_i \ge 0$ and
$\vv_i \in \Re^n$, $i=1,\ldots,r$ such that
\[
\A=\sum_{i\in [r]} \alpha_i (\vv_i)^{\otimes m} \mbox{ with }
\vv_i=\left(1,\gamma^i,\ldots,(\gamma^i)^{n-1}\right).
\]
By comparing with the $(1,\ldots,1)$-entry, we see that
$\alpha_i=0$, $i=1,\ldots,r$. This gives $\A=\OO$ which makes
a contradiction.
\ep

In the following theorem, we will show that when the order is even,
a strong Hankel tensor is indeed a limiting point of complete Hankel
tensors.

\begin{Theorem}{(\bf Completely decomposability of strong Hankel tensors)}\label{th:2}
Let $m,n \in \mathbb{N}$. Let $\A$ be an $m$th-order
$n$-dimensional strong Hankel tensor. If the order $m$ is an even
number, then $\mathcal{A}$ is a completely decomposable tensor and a
limiting point of  complete Hankel tensors. If the order $m$ is an
odd number, then $\mathcal{A}$ is
a completely $r$-decomposable tensor with
$r=(n-1)m+1$.
\end{Theorem}

\noindent {\bf Proof.}
Let $h$ be the generating function of the strong Hankel tensor $\mathcal{A}$. Then, for any $\x\in \Re^n$,
\begin{eqnarray}\label{eq:1}
f(\x):=\A \x^{\otimes m} & = &  \sum_{i_1,i_2,\ldots,i_m=1}^n v_{i_1+i_2+\ldots+i_m-m} x_{i_1}x_{i_2} \ldots x_{i_m} \nonumber \\
& = & \sum_{i_1,i_2,\ldots,i_m=1}^n  \left(\int_{-\infty}^{+\infty}  t^{i_1+i_2+\ldots+i_m-m} h(t)   dt \right)\, x_{i_1}x_{i_2} \ldots x_{i_m} \nonumber \\
&=&\int_{-\infty}^{+\infty}  \left(\sum_{i_1,i_2,\ldots,i_m=1}^n t^{i_1+i_2+\ldots+i_m-m} x_{i_1}x_{i_2} \ldots x_{i_m}  \right)h(t)dt \, \nonumber \\
&=& \int_{-\infty}^{+\infty}  \left(\sum_{i=1}^n t^{i-1} x_{i}\right)^m h(t)dt = \lim_{l \rightarrow +\infty} f_l(x),
\end{eqnarray}
where
\[
f_l(\x)=\int_{-l}^{l}  \left(\sum_{i=1}^n t^{i-1} x_{i}\right)^m h(t)dt.
\]
By the definition of Riemann integral, for each $l \ge 0$, we have
$f_l(\x)  = \lim_{k\rightarrow \infty} f_l^k(\x)$,
where $f_l^k(\x)$ is a polynomial defined by
\[
f_l^k(\x):=\sum_{j=0}^{2kl} \frac{\left(\sum_{i=1}^n
(\frac{j}{k}-l)^{i-1} x_i \right)^m h(\frac{j}{k}-l)}{k}.
\]
Fix any $l \ge 0$ and $k \in \mathbb{N}$.
Note that
\begin{eqnarray*}
f_l^k(\x):=\sum_{j=0}^{2kl} \frac{\left(\sum_{i=1}^n
(\frac{j}{k}-l)^{i-1} x_i \right)^m h(\frac{j}{k}-l)}{k}&= & \sum_{j=0}^{2kl}\left(\sum_{i=1}^n\frac{
(\frac{j}{k}-l)^{i-1}  h(\frac{j}{k}-l)^{\frac{1}{m} }}{k^{\frac{1}{m}}} x_i\right)^m \\
& = & \sum_{j=0}^{2kl} (\langle \uu_j, \x \rangle)^m,
\end{eqnarray*}
where $\uu_j=\frac{
 h\left(\frac{j}{k}-l\right)^{\frac{1}{m} }}{k^{\frac{1}{m}}}\left(1,\frac{j}{k}-l,\ldots,(\frac{j}{k}-l)^{n-1}\right)$. Here $\uu_j$ are always well-defined as
$h$ takes nonnegative values. Define $\A_l^k$ be a
symmetric tensor such that
$f_l^k(\x)=\A_l^k \x^{\otimes m}$.
Then, it is easy to see that
 each $\A_k^l$ is a complete Hankel tensor and thus a completely decomposable tensor. {  Note from
 Theorem \ref{t3.1} (iv) that} the completely decomposable tensor cone $CD_{m,n}$ is a closed convex cone when $m$ is even.
 It then follows that
$\A =\lim_{k \rightarrow \infty}\lim_{l \rightarrow
\infty}\A_k^l$ is a  {  completely decomposable tensor and a limiting point of complete Hankel tensors.}

To see the assertion in the odd order case, we use a similar
argument as in \cite{Qi15}. Pick real numbers
$\gamma_1,\ldots,\gamma_{r}$ with $r=(n-1)m+1$ and $\gamma_i \neq
\gamma_j$ for $i \neq j$. Consider the following linear equation in
$\alpha=(\alpha_1,\ldots,\alpha_r)$ with
\[
v_k=\sum_{i=1}^r \alpha_i \gamma_i^{k}, \ k=0,\ldots,(n-1)m.
\]
Note that this linear equation always has a solution say $\bar \alpha= (\bar \alpha_1,\ldots,\bar \alpha_r)$ because the matrix in the above linear equation is a
nonsingular Vandermonde matrix. Then, we see that
\[
\A_{i_1,\ldots,i_m} =v_{i_1+\ldots+i_m-m}=\sum_{i=1}^r \bar
\alpha_i \gamma_i^{i_1+\ldots+i_m-m}=\sum_{i=1}^r \bar \alpha_i
\left((\uu_i)^{\otimes m}\right)_{i_1,\ldots,i_m},
\]
where $\uu_i \in \Re^n$ is given by
$\uu_i=(1,\gamma_i,\ldots,\gamma_i^{n-1})^T$. This shows that
$\A=\sum_{i\in [r]} \bar \alpha_i (\uu_i)^{\otimes m}$.
Now, as $m$ is an odd
number, we have
\[
\A =\sum_{i=1}^r \left({\bar \alpha_i}^{\frac{1}{m}}
\uu_i\right)^{\otimes m} .
\]
Therefore, $\mathcal{A}$ is a completely decomposable tensor and the last conclusion follows.
\ep

From the preceding theorem and Proposition \ref{p4.1}, we have
the following corollary.

\begin{Corollary}{\bf (Non-closedness of the even order complete Hankel
tensor cone)} \label{example:02} When $m$ is even and $n \ge 2$, the cone which
consists of all the complete Hankel tensors of order $m$ and
dimension $n$ is not closed.   Its closure is the cone which
consists of all the strong Hankel tensors of order $m$ and dimension
$n$.
\end{Corollary}

\begin{Corollary}\label{cor:2}
Let $n \in \mathbb{N}$ and let $m$ be an even number. Let
$\mathcal{A}$ be an $m$th-order $n$-dimensional strong Hankel
tensor. Then $\A$ is an SOS tensor and a PSD tensor.
\end{Corollary}

\noindent {\bf Proof.}
This is a direct consequence of the above theorem and the fact that ${\rm CD}_{m,n} \subseteq {\rm SOS}_{m,n} \subseteq {\rm PSD}_{m,n}$ for even order $m$.
\ep

\noindent
{\bf Question 1} For Corollary \ref{example:02}, what is the
situation if the order is odd?

We have seen that the strong Hankel tensor can be regarded as a
checkable sufficient condition for positive semi-definite Hankel
tensor. We now provide a simple necessary condition for positive
semi-definite Hankel tensors.   This condition can be verified by
solving a feasibility problem of a semi-definite programming
problem. For a set $C$ with finite elements, we use $\sharp\, C$ to denote the number of elements in the set $C$.

\begin{Proposition}{\bf (SDP-type necessary condition for PSD Hankel tensors)}
Let $n \in \mathbb{N}$ and let $m$ be a positive even number. Let
$\A$ be an $m$th-order $n$-dimensional positive
semi-definite Hankel tensor defined by (\ref{e2.1}) with
generating vector $\vv = (v_0, \cdots, v_{(n-1)m})^\top$. Denote $\alpha(m,k)=\sharp\{
(i_1,\ldots,i_m): i_1+\ldots+i_m=m+k\}$, $k=0,\ldots,(n-1)m$. Then,
there exists a symmetric $\left(\frac{(n-1)m}{2}+1\right) \times
\left(\frac{(n-1)m}{2}+1\right)$ positive semi-definite matrix $Q$
such that
\begin{equation}\label{eq:SDP}
\alpha(m,k) \,  v_k = \sum_{\alpha+\beta=k} Q_{\alpha,\beta}, \ k=0,1,\ldots,(n-1)m.
\end{equation}
\end{Proposition}

\noindent {\bf Proof.}
As $\A$ is an $m$th-order $n$-dimensional positive
semi-definite Hankel tensor and $\vv$ is its
generating vector,
\[
\A\x^{\otimes m}= \sum_{i_1,i_2,\ldots,i_m=1}^n v_{i_1+i_2+\ldots+i_m-m} x_{i_1}x_{i_2} \ldots x_{i_m} \ge 0 \mbox{ for all } \x \in \Re^n.
\]
Consider $\x(t)=(x_1(t),\ldots, x_n(t))^\top$ with $x_i(t)=t^{i-1}$, $i=1,\ldots,n$ and $t \in \Re$. Then, for all $t \in \Re$, we have
\[
\phi(t)=\A \x(t)^{\otimes m}=\sum_{i_1,i_2,\ldots,i_m=1}^n v_{i_1+i_2+\ldots+i_m-m} \, t^{i_1+i_2+\ldots+i_m-m}=\sum_{k=0}^{(n-1)m} \sum_{i_1+i_2+\ldots+i_m-m=k}v_k\, t^k \ge 0.
\]
As $\phi$ is a one-dimensional polynomial which always takes
nonnegative values, $\phi$ is a sums-of-squares {  polynomial \cite{Lau,Las}.
Define $w_t=\left(1,t,\ldots,t^{\frac{(n-1)m}{2}}\right)^T$. So,}
there exists a symmetric $\left(\frac{(n-1)m}{2}+1\right) \times
\left(\frac{(n-1)m}{2}+1\right)$ positive semi-definite matrix $Q$
such that
\[
\sum_{k=0}^{(n-1)m} \alpha(m,k) v_k\, t^k=\sum_{k=0}^{(n-1)m} \sum_{i_1+i_2+\ldots+i_m-m=k}v_k\, t^k =\phi(t)= w_t^TQw_t,
\]
which is further equivalent to (by comparing the entries) the relation (\ref{eq:SDP}).
\ep

\begin{Theorem}{\bf (Positive definiteness of strong Hankel tensors)}
Let $n \in \mathbb{N}$ and let $m$ be an even number. Let
$\A$ be an $m$th-order $n$-dimensional strong Hankel
tensor. Suppose that the generating function $h$ takes positive
value almost everywhere. Then $\mathcal{A}$ is a positive definite
tensor.
\end{Theorem}

\noindent {\bf Proof.}
From (\ref{eq:1}), for any $\x\in \Re^n$,
\begin{eqnarray}
f(\x):=\A \x^{\otimes m} & = &  \int_{-\infty}^{+\infty}  \left(\sum_{i=1}^n t^{i-1} x_{i}\right)^m h(t)dt \, .
\end{eqnarray}
{  From Corollary \ref{cor:2}, $\A \x^{\otimes m} \ge 0$ for all $\x \in \Re^n$.} Suppose that there exists $\bar \x \neq \0$ such that $\A\bar \x^{\otimes m}=0$. Then,
\[
\int_{-\infty}^{+\infty}  \left(\sum_{i=1}^n t^{i-1} \bar x_{i}\right)^m h(t)dt =0.
\]
From our assumption that $h$ takes positive value almost everywhere, for each $l \ge 0$, we have
\[
\left(\sum_{i=1}^n t^{i-1} \bar x_{i}\right)^m = 0 \mbox{ for almost every } t \in [-l,l].
\]
By the continuity, this shows that
$\sum_{i=1}^n t^{i-1} \bar x_{i}  \equiv 0 \mbox{ for all } t \in [-l,l]$.
So, for each $l \ge 0$, we have
$\bar x_1+ \sum_{i=2}^n t^{i-1} \bar x_i \equiv 0 \mbox{ for all } t \in [-l,l]$.
Letting $t=0$,  we have $\bar x_1=0$.  Then,
\[
\sum_{i=2}^n t^{i-1} \bar x_i \equiv 0 \mbox{ for all } t \in [-l,l].
\]
Repeating this process, we have $\bar x_2=\ldots=\bar x_n=0$. So, $\bar \x=\0$. Therefore, we see that $\A$ is positive definite.
\ep

\begin{Theorem}{\bf (Complete positivity of strong Hankel tensors)}
Let $m,n \in \mathbb{N}$. Let $\A$ be an $m$th-order
$n$-dimensional strong Hankel tensor with a generating function $h$
on $\Re$. Suppose that $\{t \in \mathbb{R}: h(t) \neq 0\}
\subseteq \mathbb{R}_+$. Then $\A$ is a completely positive
tensor.
\end{Theorem}

\noindent {\bf Proof.}
As $\{t \in \Re : h(t) \neq 0\} \subseteq \Re_+$, we have for any $x\in \Re^n$
\begin{eqnarray}\label{eq:11}
f(\x):=\A \x^{\otimes m}
&=& \int_{-\infty}^{+\infty}  \left(\sum_{i=1}^n t^{i-1} x_{i}\right)^m h(t)dt = \int_{0}^{+\infty}  \left(\sum_{i=1}^n t^{i-1} x_{i}\right)^m h(t)dt.
\end{eqnarray}
Then, using similar line of argument as in the preceding theorem, we see that
\[
f(\x)=\lim_{l \rightarrow +\infty}\lim_{k \rightarrow \infty} f_l^k(\x),
\]
where $f_l^k(\x)$ is given by
\begin{eqnarray*}
f_l^k(\x):=\sum_{j=0}^{kl} \frac{\left(\sum_{i=1}^n
(\frac{j}{k})^{i-1} x_i \right)^m h(\frac{j}{k})}{k}&= & \sum_{j=0}^{kl}\left(\sum_{i=1}^n\frac{
(\frac{j}{k})^{i-1}  h(\frac{j}{k})^{\frac{1}{m} }}{k^{\frac{1}{m}}} x_i\right)^m
 =  \sum_{j=0}^{kl} (\langle \uu_j, \x \rangle)^m,
\end{eqnarray*}
where $\uu_j=\left(\frac{
 h(\frac{j}{k})^{\frac{1}{m} }}{k^{\frac{1}{m}}},\frac{
(\frac{j}{k})  h(\frac{j}{k}-l)^{\frac{1}{m}
}}{k^{\frac{1}{m}}},\ldots,\frac{ (\frac{j}{k})^{n-1}
h(\frac{j}{k})^{\frac{1}{m} }}{k^{\frac{1}{m}}}\right)^\top \in
\Re^n_+$ (as $h$ takes nonnegative values). Define
$\A_l^k$ be a symmetric tensor such that
$f_l^k(\x)=\A_l^k \x^{\otimes m}$. So,
 each $\A_k^l$ is a completely positive tensor. Note that the completely positive cone $CP_{m,n}$ is a closed convex cone \cite{QXX}. It then follows that
$\A=\lim_{k \rightarrow \infty}\lim_{l \rightarrow \infty}\A_k^l$ is also a  completely positive tensor.
\ep

\section{Other SOS-Hankel Tensors}
\hspace{4mm}

There are other positive semi-definite Hankel tensors, which are not strong Hankel tensors.    In \cite{Qi15}, there is such an example.   That example is for order $m=4$.  We now extend it to $m = 2k$ for any integer $k \ge 2$.  We will show that such tensors are not completely decomposable {  (and so, is also not
a strong Hankel tensor by Theorem \ref{th:2})}, but they are still SOS-Hankel tensors.

Let $m = 2k$, $n$ = 2, $k$ is an integer and $k \ge 2$.   {  Let $v_0 = v_m = 1$, $v_{2l} = v_{m-2l} = -{1 \over \left({m \atop 2l}\right)}$, $l=1,\ldots,k-1$, and $v_j = 0$ for other $j$. Let $\A = (a_{i_1\cdots i_m})$ be defined by
$a_{i_1\cdots i_m} = v_{i_1+\cdots +i_m-m},$
for $i_1, \cdots, i_m = 1, 2$.}   Then $\A$ is an even order Hankel tensor.
{  For any $\x \in \Re^2$, we have
$$\A \x^{\otimes m} = x_1^m - \sum_{j=1}^{k-1} x_1^{m-2j} x_2^{2j}   + x_2^m = \sum_{j=0}^{k-2}\left(x_1^{k-j}x_2^j - x_1^{k-j-2}x_2^{j+2}\right)^2.$$
Thus,} $\A$ is an SOS-Hankel tensor, hence a positive semi-definite Hankel tensor.  On the other hand, $\A$ is not a completely decomposable tensor.  Assume that $\A$ is a completely decomposable tensor.   Then there are vectors $\uu_j = (a_j, b_j)^\top$ for $j \in [r]$ such that
$\A = \sum_{j=1}^r \uu_j^m.$
Then for any $\x \in \Re^2$,
$$\A \x^{\otimes m} = \sum_{p=1}^r (a_px_1+b_px_2)^m = \sum_{j=0}^m \sum_{p=1}^r \left({ m \atop j}\right)a_p^{m-j}b_p^jx_1^{m-j}x_2^j.$$
{  On the other hand,
$$\A \x^{\otimes m} = x_1^m -  \sum_{j=1}^{k-1} x_1^{m-2j} x_2^{2j} + x_2^m.$$}
Comparing the coefficients of $x_1^{m-2}x_2^2$ in the above two expressions of $\A \x^m$, we have
$$\sum_{p=1}^r \left({ m \atop 2}\right)a_p^{m-2}b_p^2 = -1.$$
This is impossible.   Thus, $\A$ is not completely decomposable.


\medskip

\noindent
{\bf Question 2}  Is there an even order completely decomposable Hankel tensor, which is not a strong Hankel tensor?

\medskip

We may also construct an example for $m=6$ and $n=3$.   Let $\A \in {\rm S}_{4, 3}$ be a Hankel tensor generated by $\vv = (v_0 = \alpha, 0, \cdots, 0, v_6=1, 0, \cdots, 0, v_{12}=\alpha)^\top$. Then
\[
f(\x)\equiv \A\x^{\otimes 6} =\alpha x_1^6 + x_2^6+30x_1x_2^4x_3+90x_1^2x_2^2x_3^2+ 20x_1^3x_3^3 +\alpha x_3^6.
\]
We now show that $f$ is PSD if $\alpha \ge 480\sqrt{15}+10$.
Indeed,
\begin{equation} \label{e1}
f(\x ) = 10\left(x_1^3 +x_3^3\right)^2 + {x_2^2 \over 2} \left(x_2^2 + 30x_1x_3\right)^2 + \left[(\alpha-10)(x_1^6 +x_3^6) + {1 \over 2}x_2^6 - 360x_1^2x_2^2x_3^2 \right] \ge 0,
\end{equation}
where $\left[(\alpha-10)(x_1^6 +x_3^6) + {1 \over 2}x_2^6 - 360x_1^2x_2^2x_3^2 \right] \ge 0$ because of the arithmetic-geometric inequality.
Note that $\left[(\alpha-10)(x_1^6 +x_3^6) + {1 \over 2}x_2^6 - 360x_1^2x_2^2x_3^2 \right] \ge 0$ is a diagonal minus tail form and all positive semi-definite diagonal minus tail forms are SOS \cite[Theorem 2.3]{FK}.   Thus, if $\alpha \ge 480\sqrt{15}+10$, $f(\x)$ is also SOS.

On the other hand, we may see that the Hankel tensor $\A$ is not a strong Hankel tensor.  Let $A = (a_{ij})$ be generated by $\vv$.  Then $A$ is a $7 \times 7$ Hankel matrix, with $a_{11} = a_{77} = \alpha$, $a_{44} = a_{35} = a_{53} = a_{26} = a_{62} = a_{17} = a_{71} =1$ and $a_{ij}=0$ for other $(i, j)$.  Let $\y = (0, 0, 1, 0, -1, 0, 0, 0)^\top$.   Then $\y^\top A \y = -2 < 0$.   Hence $A$ is not PSD and $\A$ is not a strong Hankel tensor.
\bigskip

Naturally, we have the following question:

\medskip

\noindent
{\bf Question 3}  Are there PSD non-SOS-Hankel tensors?

In a certain sense, this question looks like the Hilbert problem under the Hankel constraint.

\medskip

The following question is connected with the above question.

\medskip

\noindent
{\bf Question 4}  Is the problem for determining an even order Hankel tensor is positive semi-definite or not solvable in polynomial-time?

\medskip

As discussed before, if there are no PSD non-SOS-Hankel tensors, then the problem for determining an even order Hankel tensor is positive semi-definite or not is solvable in polynomial-time.

\medskip

Before answering Question 3, we may try to answer an easier question.

\medskip

\noindent
{\bf Question 5}  Are there PSD non-SOS-Hankel tensors of order $6$ and dimension $3$?

\medskip

If there are no PSD non-SOS-Hankel tensors of order $6$ and dimension $3$, then it indicates that the Hilbert problem for PSD non-SOS polynomials has different answer if such polynomials are restricted to be Hankel polynomials.

\section{An Application}
\hspace{1mm}

\subsection{Positive Semidefinite Hankel Tensor Completion Problem}
\hspace{4mm}

An interesting problem is to fit the tensor data with prescribed Hankel structure of low rank. That is, given a tensor $\X \in {\rm S}_{m,n}$, one try to solve
the following optimization problem:
\begin{eqnarray*}
& \min_{\} \in {\rm S}_{m,n}}& \frac{1}{2}\|\A-\X \|^2 \\
&\mbox{ subject to } & \A \mbox{ is a positive semi-definite Hankel tensor} \\
& & \mathcal{A} \mbox{ is of low rank}.
\end{eqnarray*}
As a PSD Hankel tensor and the rank for tensor are hard to determine. One could consider the following alternative
\begin{eqnarray*}
 & \min_{\A \in {\rm S}_{m,n}}& \frac{1}{2}\|\A -\X \|^2 \\
&\mbox{ subject to } & \A \mbox{ is a strong Hankel tensor} \\
& & \mbox{the associated Hankel matrix of $\A$ is of low rank}.
\end{eqnarray*}
Define $l={(n-1)m+2 \over 2}$. This problem can be rewritten as
\begin{eqnarray*}
& \min_{v \in \mathbb{R}^{(n-1)m+1}}& \frac{1}{2}\sum_{i_1,\ldots,i_m=1}^n\|v_{i_1+\ldots+i_m-m}-\X_{i_1,\ldots,i_m}\|^2 \\
&\mbox{ subject to } & A_{\alpha,\beta}=v_{\alpha+\beta-2}, \ \alpha,\beta=1,\ldots,l, \\
& & A \mbox{ is a positive semi-definite matrix} \mbox{ and } A
\mbox{ is of low rank}.
\end{eqnarray*}
The low rank constraint is often nonconvex and nonsmooth, and so, the problem is still a hard problem to solve. As the trace norm promotes a low rank solution, one popular approach is to solve its nuclear norm heuristic
(see, for example \cite{JMZ})  the following form:
\begin{eqnarray*}
& \min_{v \in \mathbb{R}^{(n-1)m+1}}& \frac{1}{2}\sum_{i_1,\ldots,i_m=1}^n\|v_{i_1+\ldots+i_m-m}-\mathcal{X}_{i_1,\ldots,i_m}\|^2 + \mu \|A\|_{\rm tr} \\
&\mbox{ subject to } & A_{\alpha,\beta}=v_{\alpha+\beta-2}, \ \alpha,\beta=1,\ldots,l, \mbox{ and } A \succeq 0,
\end{eqnarray*}
where $\|A\|_{\rm tr}$ denotes the  trace norm of $A$ and is defined as the sum of all eigenvalues of $A$. Define two linear maps $M:\Re^{(n-1)m+1} \rightarrow {\rm S}^l:={\rm S}_{2,l}$ and $\PP :\Re^{(n-1)m+1} \rightarrow {\rm S}_{m,n}$ by
\[
M \, v=\big(v_{\alpha+\beta-2}\big)_{1\le \alpha,\beta \le r} \mbox{ and } (\PP(\vv)_{i_1,\ldots,i_m}= v_{i_1+\ldots+i_m-m}.
\]
Then, the trace norm  problem can be further simplified as
\begin{eqnarray*}
(TCP) & \displaystyle \min_{\vv \in \Re^{(n-1)m+1}, A \in S^l }& \frac{1}{2} \|\PP(\vv)-\X \|^2 + \mu \|A\|_{\rm tr} \\
&\mbox{ subject to } & A-M\, \vv=0, \mbox{ and } A \succeq 0.
\end{eqnarray*}
The associated augmented Lagrangian for $(TCP)$ can be formulated as
\[
L_{\rho}(\vv,A,Z)=\frac{1}{2} \|\PP(\vv)-\X \|^2 +  \mu \|A\|_{\rm tr}+ {\rm Tr}[Z (A-M\, \vv)]+\frac{\rho}{2}\|A-M \vv\|^2.
\]
We now propose an alternating direction method of multiplier (ADMM) to solve $(TCP)$.

\bigskip

{\bf ADMM for solving $(TCP)$}

Step 0. Given $\mathcal{X} \in S_{m,n}$. Choose an initial point $(v^0,A^0,Z^0) \in \mathbb{R}^{(n-1)m+1} \times S^l \times S^l$ and $\rho>0$. Set $k=1$.

Step 1. Find $\vv^k={\rm argmin}_\vv\{L_{\rho}(\vv,A^{k-1},Z^{k-1})\}$.

Step 2. Find $A^k={\rm argmin}_A\{L_{\rho}(v^k,A,Z^{k-1}): A \succeq 0\}$.

Step 3. Let $Z^k=Z^{k-1}+\rho (A^k -M \, v^k)$.

Step 4. Let $k=k+1$ and go back to Step 1.

\bigskip

The computational cost of (ADMM) is not heavy which makes it suitable for solving large size problem. In fact, we note that, in Step 1, $\vv^k$ can be found by solving a linear equation in $v$:
\[
(\PP^*\PP+\rho M^*M) v = \PP^*\X +M^*(\rho A^{k-1}+Z^{k-1}),
\]
where $\PP^*$ and $M^*$ are the adjoint mappings of $\PP$ and $M$ respectively.
Moreover, in Step 2, $A^k$ indeed has a closed form solution as
\begin{eqnarray*}
A^k &= & {\rm argmin}_A\{L_{\rho}(\vv^k,A,Z^{k-1}): A \succeq 0\} \\
&= & {\rm argmin}_A\left\{\mu \|A\|_{\rm tr}+ {\rm Tr}[Z^{k-1} (A-M\, \vv)]+\frac{\rho}{2}\|A-M \vv^k\|^2: A \succeq 0\right\} \\
& = & {\rm argmin}_A\left\{\mu {\rm Tr}[I_l A] + {\rm Tr}[Z^{k-1} A]+\frac{\rho}{2}\|A-M \vv^k\|^2: A \succeq 0\right\} \\
& = & {\rm argmin}_A\left\{\frac{\rho}{2}\|A+ \frac{1}{\rho}\left(\mu I_l+Z^{k-1}\right)-M \vv^k\|^2: A \succeq 0 \right\} \\
& = & \sum_{i=1}^l \max\{\sigma_i,0\} \uu_i\uu_i^\top
\end{eqnarray*}
where $\sigma_i$ and $\uu_i$ are obtained from the SVD decomposition of $M \,v^k-\frac{1}{\rho}(\mu I_l+Z^{k-1})$, that is,
\[
M \, \vv^k-\frac{1}{\rho}\left(\mu I_l+Z^{k-1}\right)= \sum_{i=1}^l \sigma_i \uu_i\uu_i^\top.
\]
The convergence of the ADMM method has been well-studied by a lot of researchers. For simplicity purpose, we omit the details and refer the interested reader to  \cite{HY}.
\subsection{Numerical Tests}
\hspace{4mm}

{  To illustrate the algorithm of ADMM for solving (TCP), we first generate two random instances of symmetric tensors. Then, we solve the
corresponding positive semidefinite Hankel tensor completion problem via the proposed ADMM algorithm. In our numerical test, we set $\mu$ and $\rho$ in the
ADMM algorithm as $0.1$ and $10$ respectively.  Our proposed
algorithm works very well by using these parameters.

\begin{example}
Consider a $4$th-order $3$-dimensional  symmetric tensor $\mathcal{X}$ given by

$\mathcal{X}(:,:,1,1) =\left(
               \begin{array}{ccc}
   -0.2972  &  0.4307  &  0.4444\\
    0.4307  & -0.4029 &  -0.0274\\
    0.4444   &-0.0274  &  0.0647\\
               \end{array}
             \right),$
$\mathcal{X}(:,:,2,1) =\left(
               \begin{array}{ccc}
    0.4307&   -0.4029 &  -0.0274\\
   -0.4029  &  0.1085  &  0.1760\\
   -0.0274  &  0.1760  & -0.2574\\
              \end{array}
             \right),$

$\mathcal{X}(:,:,3,1) =\left(
               \begin{array}{ccc}
    0.4444  & -0.0274  &  0.0647\\
   -0.0274 &   0.1760  & -0.2574\\
    0.0647  & -0.2574 &  -0.3208\\
              \end{array}
             \right),$
$\mathcal{X}(:,:,1,2) =\left(
               \begin{array}{ccc}
    0.4307  & -0.4029  & -0.0274\\
   -0.4029  &  0.1085 &   0.1760\\
   -0.0274  &  0.1760 &  -0.2574\\
              \end{array}
             \right),$

$\mathcal{X}(:,:,2,2) =\left(
               \begin{array}{ccc}
   -0.4029  &  0.1085 &   0.1760\\
    0.1085  &  0.9152  & -0.0821\\
    0.1760  & -0.0821 &  -0.2815\\

              \end{array}
             \right),$
$\mathcal{X}(:,:,3,2) =\left(
               \begin{array}{ccc}
   -0.0274 &   0.1760 &  -0.2574\\
    0.1760  & -0.0821 &  -0.2815\\
   -0.2574 &  -0.2815 &   0.2773\\
              \end{array}
             \right),$

$\mathcal{X}(:,:,1,3) =\left(
               \begin{array}{ccc}
    0.4444 &  -0.0274  &  0.0647\\
   -0.0274  &  0.1760 &  -0.2574\\
    0.0647  & -0.2574 &  -0.3208\\
              \end{array}
             \right),$
$\mathcal{X}(:,:,2,3) =\left(
               \begin{array}{ccc}
   -0.0274  &  0.1760  & -0.2574\\
    0.1760 &  -0.0821  & -0.2815\\
   -0.2574 &  -0.2815  &  0.2773\\
              \end{array}
             \right),$

$\mathcal{X}(:,:,3,3) =\left(
               \begin{array}{ccc}
    0.0647  & -0.2574  & -0.3208\\
   -0.2574 &  -0.2815  &  0.2773\\
   -0.3208 &   0.2773  & -0.5347\\
              \end{array}
             \right).$

Solving the strong Hankle tensor completion problem for $\mathcal{X}$ via the proposed ADMM method,  we obtain $$v=(0.0086,    0.0056,    0.0036,    0.0022,    0.0014,    0.0009,    0.0006,    0.0004,    0.0002)^T$$ and the associated Hankel matrix  $A= (A_{ij})_{1 \le i \le j \le 5}$ with $A_{ij}=v_{i+j-2}$
is of rank $1$.
\end{example}

\begin{example}
Consider a given $4$th-order $3$-dimensional  symmetric tensor $\mathcal{X}$ given by

$\mathcal{X}(:,:,1,1) =\left(
               \begin{array}{ccc}
   -0.7384 &   0.2309  &  0.3538\\
    0.2309  & -0.4025 &   0.2401\\
    0.3538 &   0.2401 &  -0.2167\\
               \end{array}
             \right),$
$\mathcal{X}(:,:,2,1) =\left(
               \begin{array}{ccc}
    0.2309 &  -0.4025  &  0.2401\\
   -0.4025 &   0.1324  & -0.1888\\
    0.2401 &  -0.1888  & -0.1495\\
             \end{array}
             \right),$

$\mathcal{X}(:,:,3,1) =\left(
               \begin{array}{ccc}
    0.3538  &  0.2401  & -0.2167\\
    0.2401  & -0.1888 &  -0.1495\\
   -0.2167  & -0.1495  &  0.3234\\
              \end{array}
             \right),$
$\mathcal{X}(:,:,1,2) =\left(
               \begin{array}{ccc}
    0.2309  & -0.4025  &  0.2401\\
   -0.4025  &  0.1324 &  -0.1888\\
    0.2401  & -0.1888 &  -0.1495\\
              \end{array}
             \right),$

$\mathcal{X}(:,:,2,2) =\left(
               \begin{array}{ccc}
   -0.4025  &  0.1324  & -0.1888\\
    0.1324  & -0.3712  &  0.0019\\
   -0.1888  &  0.0019  & -0.1546\\

              \end{array}
             \right),$
$\mathcal{X}(:,:,3,2) =\left(
               \begin{array}{ccc}
    0.2401 &  -0.1888  & -0.1495\\
   -0.1888  &  0.0019  & -0.1546\\
   -0.1495  & -0.1546 &  -0.0395\\
              \end{array}
             \right),$

$\mathcal{X}(:,:,1,3) =\left(
               \begin{array}{ccc}
    0.3538 &   0.2401  & -0.2167\\
    0.2401  & -0.1888  & -0.1495\\
   -0.2167 &  -0.1495  &  0.3234\\
              \end{array}
             \right),$
$\mathcal{X}(:,:,2,3) =\left(
               \begin{array}{ccc}
    0.2401  & -0.1888  & -0.1495\\
   -0.1888  &  0.0019 &  -0.1546\\
   -0.1495 &  -0.1546 &  -0.0395\\
              \end{array}
             \right),$

$\mathcal{X}(:,:,3,3) =\left(
               \begin{array}{ccc}
   -0.2167  & -0.1495   & 0.3234\\
   -0.1495  & -0.1546   &-0.0395\\
    0.3234  & -0.0395  &  0.9162\\
              \end{array}
             \right).$

Solving the strong Hankle tensor completion problem for $\mathcal{X}$ via the proposed ADMM method,  we obtain
$$v=( 0, 0, -0.0001, 0.0003, 0.0001, -0.0024, 0.0120, -0.0390, 0.7741)^T$$ and and the associated Hankel matrix  $A= (A_{ij})_{1 \le i \le j \le 5}$ with $A_{ij}=v_{i+j-2}$
is of rank $2$.
\end{example}

These preliminary numerical results show that  the algorithm of ADMM for solving (TCP) is efficient.}

\hspace{4mm}


\end{document}